\ifcase2
\or
\documentclass[12pt,a5paper]{article}
\def\atPreamble{\usepackage{mybiblatex}}
\or
\documentclass[12pt,a4paper]{article}
\AtEndDocument{\bibliography{ajr,bib}}
\def\atPreamble{%
    \RaisedNamesfalse
    \usepackage{natbib,twoopt}
    \bibliographystyle{agsm}
    \let\harvardurl\url
    \newcommandtwoopt{\footcite}[3][][]%
        {\footnote{\cite[##1][##2]{##3}}}
    }
\fi

\title{Embed to rigorously and accurately homogenise quasi-periodic multi-scale heterogeneous PDEs, with computer algebra}
\author{A.~J. Roberts
\thanks{University of Adelaide, South Australia. 
\protect\url{http://orcid.org/0000-0001-8930-1552},
\protect\url{mailto:profajroberts@protonmail.com}}
}
\date{\today}

\IfFileExists{ajr.sty}{\usepackage{ajr}}{}
\usepackage[defaultlines=2,all]{nowidow}
\usepackage{pgfplots}
\pgfplotsset{compat=newest}
\usepgfplotslibrary{external}
\usepackage{amsmath,amssymb,defns,reducecode,mycleveref}

\atPreamble

\renewcommand{\vec}[1]{\text{\boldmath$#1$}}
\def\Res{\operatorname{Res}}
\def\rve{\textsc{rve}}
\def\cas{\textsc{cas}}
\def\ov{\vec 0}
\Vec{phi}
\Cal C\Cal D\Cal K\Cal L
\Bb E\Bb H\Bb X
\Frak D\Frak K\Frak L\Frak p\Frak u\Frak U\Frak v
\allowdisplaybreaks

\begin{document}

\maketitle

\begin{abstract}
For microscale heterogeneous \pde{}s, this article further develops novel theory and methodology for their macroscale mathematical\slash asymptotic homogenization.
This article specifically encompasses the case of quasi-periodic heterogeneity with finite scale separation: no scale separation limit is required.
Dynamical systems theory frames the homogenization as a slow manifold of the ensemble of all phase-shifts of the heterogeneity.
Depending upon any perceived scale separation \emph{within} the quasi-periodic heterogeneity, the homogenization may be done in either one step, or two sequential steps: the results are equivalent.
The theory not only assures us of the existence and emergence of the homogenization, it also provides a practical systematic method to construct the homogenization to any specified order.
For a class of heterogeneities, we show that the macroscale homogenization is potentially valid down to lengths which are just twice that of the microscale heterogeneity!
This methodology provides a new rigorous and flexible approach to homogenization that potentially also provides correct initial and boundary conditions, treatment of forcing and control, and analysis of~uncertainty.
\end{abstract}

\tableofcontents

\section{Introduction}

Many engineering and scientific structures have microscale structure, such as lattice materials \cite[e.g.,][]{Somnic2022, Somnic2022b, Barchiesi2020}, the windings in electrical machinery \cite[e.g.,][]{Romanazzi2016}, 
electromagnetism in micro-structured material \cite[e.g.,][]{Craster2015, Niyonzima2014},  food drying \cite[e.g.,][]{Welsh2018}, and
slow Stokes flow through porous media \cite[e.g.,][]{Brown2011}.
The engineering challenge is to understand the dynamics on a scale significantly larger than the micro-structure.
Homogenization, via the recognitionn of multiple physical length scales, is the common approach\footcite[e.g.,][]{Gustafsson03, Engquist08, Cooper2018, Abdulle2020c}.
This article further develops a new approach to modelling the large scale homogenized dynamics via a newly established general rigorous theory \citep[\S2.5]{Roberts2013a}.
This new approach provides a route to systematic refinements of the homogenization, without any of the usual assumptions except a \emph{finite} separation of scales, and to novel quantification of the remainder error.

Consider materials with complicated microstructure:  
we want to model their large scale dynamics by equations with effective, `average' coefficients.
Heterogeneous diffusion in 1D is the simplest such example: suppose the material, in spatial domain~\([0,L]\), has structure so that `material' or `heat'~\(u(t,x)\) diffuses according to 
\begin{equation}
\D tu=\D x{}\left\{\kappa(x)\D xu\right\},
\quad 0<x<L\,,
\label{eq:hdifpde}
\end{equation}
where the heterogeneous diffusion coefficient~\(\kappa(x)\) has some complicated multiscale structure.
For definiteness, let's initially suppose the macroscale boundary conditions are simply that~\(u\) is \(L\)-periodic in space~\(x\).
Our challenge is to derive the model \pde
\begin{equation}
\D tU=K\DD xU+\cdots
\label{eq:edifpde}
\end{equation}
for some effective mean field~\(U(x,t)\), some effective macroscale diffusivity~\(K\), with some higher order corrections~(\(\cdots\)), and possibly an error estimate~(\(\cdots\)).
We also comment throughout on the closely related problem\footcite[e.g.,][]{Shahraki2020, Cornaggia2020, Abdulle2020c} of wave propagation through heterogeneous material, \(u_{tt}=\partial_x\{\kappa(x)u_x\}\), and its corresponding effective homogenized \pde\ models \(U_{tt}=KU_{xx}+\cdots\) that predict its macroscale dispersive waves.

Here we focus on cases where the heterogeneity~\(\kappa(x)\) is quasi-periodic in space.
For example, \cref{secxchm,secsrh} focuses on homogenization for the family of heterogeneous diffusivities \(\kappa=1/[1+a_1\cos(2\pi x/\ell_1)+a_2\cos(2\pi x/\ell_2)]\) for arbitrary microscale periods~\(\ell_1,\ell_2\) with \(\ell_1>\ell_2\).
\cref{secpse,secsmh} develops theory and supporting characteristics of  homogenization in the general case where the heterogeneity~\(\kappa\) is expressible as a function of two variables, say~\(\kappa(x_1,x_2)\), in which~\(\kappa\) is \(\ell_i\)-periodic in variable~\(x_i\). 
Then the quasi-periodic heterogeneity in~\eqref{eq:hdifpde} is~\(\kappa(x,x)\).

Almost all extant conventional approaches to homogenization depend upon the identification of a \emph{representative volume element}~(\rve) \cite[e.g.,][]{Somnic2022b}, or a representative cell.
Even in the simple case of periodic media an uncertain issue is an appropriate size of an \rve\footcite[e.g.,][\S2.3, \S1.2.2, \S5.1, p.2, p.2, p.13, resp.]{VinhPhuNguyen2011, Saeb2016, Matous2017, Klarmann2019, Schneider2021, Somnic2022}.
An issue with such \rve-like approaches is that for quasi-periodic heterogeneities there is \emph{no} one \rve\ \cite[e.g.,][\S2.1]{Anthoine2010}.  
When the microscale lengths~\(\ell_1,\ell_2\) are irrationally related, then all volume elements, no matter what size is chosen, are different: no volume element is representative.
Relatedly, when the microscale lengths~\(\ell_1,\ell_2\) have a large least common multiple, then any true \rve\ must be large---unnecessarily large.
Instead \cref{secpse,secsmh} develop a novel a rigorous approach that, in effect, uses an `\rve' which is the size of the largest microscale length, albeit of more spatial dimensions.

I expect straightforward extensions of the approach developed herein to empower cognate analysis and modelling of the following:
\begin{itemize}
\item cases where the heterogeneity is quasi-periodic in any number of microscale lengths~\(\ell_1>\ell_2>\ell_3>\cdots\)\,;
\item \pde{}s~\eqref{eq:hdifpde} where the field~\(u\) is itself in a vector space such as the cross-sectional structure of either shear dispersion \cite[e.g.,][]{Mercer90, Mercer94a}, or of elastic beams \cite[e.g.,][]{Roberts93};
\item multiple macroscale space dimensions as in the modelling of plates and shells \cite[e.g.,][]{Roberts2016a};
\item \pde{}s like~\eqref{eq:hdifpde} but with macroscale variations in the material parameters as well as microscale variations \cite[e.g.,][]{Anthoine2010}---sometimes called \emph{functionally graded materials}; 
\item \pde{}s like~\eqref{eq:hdifpde} but with nonlinear effects \cite[e.g.,][]{Matous2017}. 
\end{itemize}

Here we significantly extend the new approach to modelling the multiscale class of systems where notionally both microscales \(\ell_1,\ell_2\ll L\), the macroscale length.
But the \cref{sec3sh} shows that the analysis and results also include three-scale multiscale systems where notionally \(\ell_2\ll\ell_1\ll L\), that is, where there is a microscale, a mesoscale, and a macroscale\footcite[e.g.,][]{RamirezTorres2018, CruzGonzalez2022}.

Analytic techniques such as mathematical homogenization or the method of multiple scales\footcite[e.g.,][]{Shahraki2020, Cornaggia2020} provide a mechanism for deriving effective \pde\ models such as~\eqref{eq:edifpde}.
But such previous underlying theory requires the limit of infinite scale separation between the microscale `cell'-sizes~\(\ell_i\) and the macroscale domain size~\(L\)\footcite[e.g.,][p.5, \S2.4 \& p.6, resp.]{Welsh2018, Barchiesi2020, Altmann2021}.
In contrast, and arising from developments in dynamical systems, \cref{secpse,secsmh} respectively use phase-shift embedding and slowly varying theory \cite[]{Roberts2013a, Roberts2016a} to provide rigorous support to the modelling at finite scale separation: that is, when the macroscale length scale~\(L\) is larger than the microscales~\(\ell_1\) and~\(\ell_2\), but not infinitely so.
In the approach here, not only is there no limit, there is no need for a defined~``\(\epsilon\)''.
Such finite scale separation empowers this approach to capture the dispersion in wave homogenization which \cite{Abdulle2020c} call for in their discussion~[p.3] that ``new effective models are required''.
The phase-shift embedding developed here further illustrates the general principle that macroscale models of microscale structures are often best phrased \text{as ensemble average}s\footcite[e.g.,][]{vanKampen92, Drugan96, Birch2006}.

It is well recognised that homogenization is physically effective for reasonably well separated length scales: that the resolved macroscale structures, with length~\(\ell_M\) say, satisfy \(\ell_M\gg\ell_1\).
But what does this mean?  
Much extant homogenization theory requires the mathematical limit \(\ell_1/\ell_M\to0\).
However, \emph{physically} the ratio~\(\ell_1/\ell_M\) is always finite.
By setting homogenization on a new rigorous footing, \cref{secsrh} analyses homogenization to high-order and presents evidence that, in a class of problems, a \emph{quantitative} limit to the spatial resolution of homogenization is that \(\ell_1/\ell_M\lesssim1/2\), equivalently \(\ell_M\gtrsim 2\ell_1\).
Mathematically, a wider scale separation is \text{not necessary for homogenization.}

Of course, practically we usually prefer a wider scale separation because of many complicating factors in construction and use of practical homogenization models.
But, in the sense discussed by \cref{secsrh}, this quantitative lower bound on the required scale separation indicates how little scale separation is necessary.
Indeed, because it applies at finite scale separation, the analysis here encompasses \emph{high-contrast media} and so also detects features occurring in this particularly interesting class of problems \cite[e.g.,][p.2]{Cooper2018}.

The new approach to homogenization that we develop here creates a powerful dynamical systems framework for addressing other modelling issues in the future.
For a time-dependent simulation of homogenized dynamics one must provide some initial conditions for the \pde~\eqref{eq:edifpde}.  
These initial conditions are surprisingly nontrivial at fininte scale separationn, but the new framework developed here comes with well-established methods to derive the correct initial conditions for long-term forecasts\footcite[e.g.,][]{Roberts89b, Roberts93, Cox93b, Roberts97b}.
The projection that provides the correct initial conditions, also provides the correct projection of any applied `forcing' \cite[\S7]{Roberts89b}, including the cases of control, and system uncertainty. 
System uncertainty is a challenge for other techniques\footcite[e.g.,][\S8 \& p.3, resp.]{Matous2017, Altmann2021}.
Further, the homogenized \pde~\eqref{eq:edifpde} needs boundary conditions---a ``highly challenging problem'' according to \cite{Shahraki2020}, and which \cite{Cornaggia2020} called a ``complex problem that is still an active research topic''.  
Such boundary conditions are nontrivial, but again the new framework developed here comes with a rigorous approach to derive correct boundary conditions for accurate predictions at finite scale separation\footcite[e.g.,][]{Roberts92c, Roberts93, Mercer94a, Chen2014, Chen2015, Chen2016}.
That is, this approach to homogenization potentially provides a rational complete modelling of not just the \pde\ at finite scale separation, but also material variations, uncertainty, control, forcing, and initial and boundary conditions.
In their review of nonlinear homogenization, \cite{Matous2017} [\S8] commented that ``verification ensures both order of accuracy and consistency \ldots\ even mathematical foundations for verification are lacking''---the approach developed here provides the required mathematical foundation.

The approach herein has no need of a variational formulation of the governing equations\footcite[e.g.,][\S2.4 \& \S3, resp.]{Barchiesi2020, Altmann2021}, and so applies to a much wider class of problems than many homogenization methods.

Throughout this article, in the details of the construction of homogenizations, one may observe that the variables~\(x_0,x_1,x_2\) appear in much the same sort of algebraic expressions traditionally obtained by the method of multiple scales.
What this parallelism implies is that in suitable circumstances the algebraic machinations of the method of multiple scales are indeed sound.
What is new here is that we place such algebraic details in a new and more rigorous framework.
The new framework developed here provides clarity, precision, power, and greater physical interpretability to homogenization.

\section{Phase-shift embedding}
\label{secpse}

Let's embed the specific given \pde~\eqref{eq:hdifpde} into a family of \pde\ problems formed by all phase-shifts of the quasi-periodic microscale.
This embedding is a novel and rigorous twist to the concept of a \emph{representative volume element}.

\begin{figure}
\centering
\caption{\label{fig:pattensemble}%
schematic cylindrical domain~\cD\ of the multiscale embedding \pde~\eqref{eq:emdifpde} for a field \(\fu(t,x_0,x_1,x_2)\) (here \(L=\ell_0=6\), \(\ell_1=1.62\), \(\ell_2=0.72\)).  
We obtain solutions of the heterogeneous diffusion \pde~\eqref{eq:hdifpde}, or~\eqref{eq:shdifpde}, on such blue lines as \(u(t,x)=\fu(t,x,x+\phi_1,x+\phi_2)\) for every pair of constant phases~\(\phi_1\) and~\(\phi_2\) (here \(\phi_1=0.82\) and \(\phi_2=0.32\)), and where the last two arguments of~\fu\ are modulo~\(\ell_1\) and~\(\ell_2\), respectively.}
%
\tikzsetnextfilename{Figs/phaseShiftEmbedding}
\begin{tikzpicture}

\begin{axis}[%
xmin=0,
xmax=6,
tick align=outside,
xlabel style={font=\color{white!15!black}},
xlabel={$x_0,x$},
ymin=0.00186145504724156,
ymax=1.61792943254703,
ylabel style={font=\color{white!15!black}},
ylabel={$x_1$},
zmin=0.000318087232761055,
zmax=0.722745278232677,
zlabel style={font=\color{white!15!black}},
zlabel={$x_2$},
view={-24}{23},
axis background/.style={fill=white},
axis equal image=true,width=120mm
]
\addplot3 [color=blue, mark=o, mark options={solid, blue}]
 table[row sep=crcr] {%
0	0.821997410046821	0.315958873732635\\
0.402	1.22399741004682	0.717958873732635\\
};
 \addplot3 [color=blue, mark=o, mark options={solid, blue}]
 table[row sep=crcr] {%
0.408	1.22999741004682	0.000352075982655631\\
0.792	1.61399741004682	0.384352075982656\\
};
 \addplot3 [color=blue, mark=o, mark options={solid, blue}]
 table[row sep=crcr] {%
0.798	0.00196342129692595	0.390352075982656\\
1.128	0.331963421296926	0.720352075982656\\
};
 \addplot3 [color=blue, mark=o, mark options={solid, blue}]
 table[row sep=crcr] {%
1.134	0.337963421296926	0.00274527823267667\\
1.854	1.05796342129693	0.722745278232677\\
};
 \addplot3 [color=blue, mark=o, mark options={solid, blue}]
 table[row sep=crcr] {%
1.86	1.06396342129693	0.00513848048269794\\
2.412	1.61596342129693	0.557138480482698\\
};
 \addplot3 [color=blue, mark=o, mark options={solid, blue}]
 table[row sep=crcr] {%
2.418	0.00392943254703138	0.563138480482698\\
2.574	0.159929432547031	0.719138480482697\\
};
 \addplot3 [color=blue, mark=o, mark options={solid, blue}]
 table[row sep=crcr] {%
2.58	0.165929432547031	0.00153168273271875\\
3.3	0.885929432547031	0.721531682732719\\
};
 \addplot3 [color=blue, mark=o, mark options={solid, blue}]
 table[row sep=crcr] {%
3.306	0.891929432547031	0.0039248849827398\\
4.02	1.60592943254703	0.71792488498274\\
};
 \addplot3 [color=blue, mark=o, mark options={solid, blue}]
 table[row sep=crcr] {%
4.026	1.61192943254703	0.000318087232761055\\
4.032	1.61792943254703	0.00631808723276128\\
};
 \addplot3 [color=blue, mark=o, mark options={solid, blue}]
 table[row sep=crcr] {%
4.038	0.00589544379713658	0.0123180872327615\\
4.746	0.713895443797137	0.720318087232762\\
};
 \addplot3 [color=blue, mark=o, mark options={solid, blue}]
 table[row sep=crcr] {%
4.752	0.719895443797136	0.0027112894827821\\
5.472	1.43989544379714	0.722711289482783\\
};
 \addplot3 [color=blue, mark=o, mark options={solid, blue}]
 table[row sep=crcr] {%
5.478	1.44589544379714	0.00510449173280314\\
5.646	1.61389544379714	0.173104491732803\\
};
 \addplot3 [color=blue, mark=o, mark options={solid, blue}]
 table[row sep=crcr] {%
5.652	0.00186145504724156	0.179104491732804\\
6	0.349861455047241	0.527104491732803\\
};
 \addplot3 [color=red, dashed]
 table[row sep=crcr] {%
0.402	1.22399741004682	0.717958873732635\\
0.408	1.22999741004682	0.000352075982655631\\
};
 \addplot3 [color=red, dashed]
 table[row sep=crcr] {%
0.792	1.61399741004682	0.384352075982656\\
0.798	0.00196342129692595	0.390352075982656\\
};
 \addplot3 [color=red, dashed]
 table[row sep=crcr] {%
1.128	0.331963421296926	0.720352075982656\\
1.134	0.337963421296926	0.00274527823267667\\
};
 \addplot3 [color=red, dashed]
 table[row sep=crcr] {%
1.854	1.05796342129693	0.722745278232677\\
1.86	1.06396342129693	0.00513848048269794\\
};
 \addplot3 [color=red, dashed]
 table[row sep=crcr] {%
2.412	1.61596342129693	0.557138480482698\\
2.418	0.00392943254703138	0.563138480482698\\
};
 \addplot3 [color=red, dashed]
 table[row sep=crcr] {%
2.574	0.159929432547031	0.719138480482697\\
2.58	0.165929432547031	0.00153168273271875\\
};
 \addplot3 [color=red, dashed]
 table[row sep=crcr] {%
3.3	0.885929432547031	0.721531682732719\\
3.306	0.891929432547031	0.0039248849827398\\
};
 \addplot3 [color=red, dashed]
 table[row sep=crcr] {%
4.02	1.60592943254703	0.71792488498274\\
4.026	1.61192943254703	0.000318087232761055\\
};
 \addplot3 [color=red, dashed]
 table[row sep=crcr] {%
4.032	1.61792943254703	0.00631808723276128\\
4.038	0.00589544379713658	0.0123180872327615\\
};
 \addplot3 [color=red, dashed]
 table[row sep=crcr] {%
4.746	0.713895443797137	0.720318087232762\\
4.752	0.719895443797136	0.0027112894827821\\
};
 \addplot3 [color=red, dashed]
 table[row sep=crcr] {%
5.472	1.43989544379714	0.722711289482783\\
5.478	1.44589544379714	0.00510449173280314\\
};
 \addplot3 [color=red, dashed]
 table[row sep=crcr] {%
5.646	1.61389544379714	0.173104491732803\\
5.652	0.00186145504724156	0.179104491732804\\
};
 \end{axis}
\end{tikzpicture}%
\end{figure}
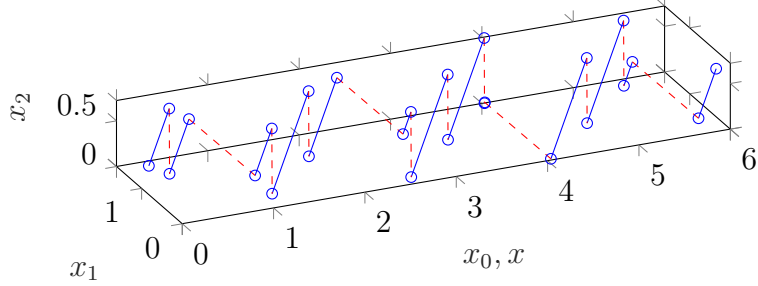%
As indicated by \cref{fig:pattensemble}, let's create the desired embedding by considering a field~\(\fu(t,x_0,x_1,x_2)\) satisfying the \pde
\begin{equation}
\D t\fu=\left(\D {x_0}{}+\D {x_1}{}+\D {x_2}{}\right)\left\{\kappa({x_1},x_2)\left(\D {x_0}\fu+\D {x_1}\fu+\D {x_2}\fu \right)\right\}
\label{eq:emdifpde}
\end{equation}
in the domain \(\cD:=\{(x_0,x_1,x_2): {x_i}\in[0,\ell_i]\}\), for the macroscale \(\ell_0=L\) and microscale~\(\ell_1,\ell_2\), and with \(\ell_i\)-periodic boundary conditions in~\({x_i}\).
We define \(\kappa_{\min}:=\min_{{x_i}\in[0,\ell_i]}\kappa(x_1,x_2)\),
and assume \(\kappa_{\min}>0\).
The domain~\cD\ in \(x_0x_1x_2\)-space (\cref{fig:pattensemble}) is often termed \emph{cylindrical} as it is long and thin.
I emphasise that this domain has finite aspect ratio: we do \emph{not} take any limit involving an aspect ratio tending to zero nor to infinity.

\cref{fig:pattensemble} indicates that we consider \(x_0=x\)\,.
The distinction between \(x_0\)~and~\(x\) is that partial derivatives in~\(x_0\) are done keeping the other~\(x_i\) constant, whereas partial derivatives in~\(x\) are done keeping the phases~\(\phi_i\) constant.

\begin{lemma}\label{lemeqv}
For every solution~\(\fu(t,x_0,x_1,x_2)\) of the embedding \pde~\eqref{eq:emdifpde}, and for every vector of phases \(\phiv:=(\phi_1,\phi_2)\), the field \(u_\phiv(t,x):=\fu(t,x,x+\phi_1,x+\phi_2)\) (that is, the field~\(\fu\) evaluated on the solid-blue lines in \cref{fig:pattensemble}) satisfies the heterogeneous diffusion \pde
\begin{equation}
\D t{u_\phiv}=\D x{}\left\{ \kappa(x+\phi_1,x+\phi_2)\D x{u_\phiv} \right\}.
\label{eq:shdifpde}
\end{equation}
Hence the field \(u(t,x):=\fu(t,x,x,x)\) satisfies the given heterogeneous \pde~\eqref{eq:hdifpde}.
\end{lemma}

Recall that the most common boundary conditions assumed for \rve{}s are periodic, whereas here the \(\ell_1,\ell_2\)-periodic boundary conditions are not assumed but arise naturally due to the ensemble of phase-shifts.
That is, what previously had to be assumed, here arises naturally.

\begin{proof}
Start by considering the left-hand side time derivative in \pde~\eqref{eq:shdifpde}:
\begin{align*}
\D t{u_\phiv}&
=\D t{}\fu(t,x,x+\phi_1,x+\phi_2)
\\&
=\left[\D t\fu\right]_{(t,x,x+\phi_1,x+\phi_2)}
\quad\text{(which by \pde~\eqref{eq:emdifpde} becomes)}
\\&
=\left[\left(\D{x_0}{}+\D {x_1}{}+\D {x_2}{}\right)\left\{\kappa({x_1},x_2)\left(\D{x_0}\fu+\D {x_1}\fu+\D {x_2}\fu \right)\right\}\right]_{(t,x,x+\phi_1,x+\phi_2)}
\\&
=\D x{}\left\{\left[ \kappa({x_1},x_2)\left(\D{x_0}\fu+\D {x_1}\fu+\D {x_2}\fu \right)\right]_{(t,x,x+\phi_1,x+\phi_2)}\right\}
\\&
=\D x{}\left\{\kappa({x+\phi_1},x+\phi_2)\left[\D{x_0}\fu+\D {x_1}\fu+\D {x_2}\fu \right]_{(t,x,x+\phi_1,x+\phi_2)}\right\}
\\&
=\D x{}\left\{\kappa({x+\phi_1},x+\phi_2)\D x{}\fu(t,x,x+\phi_1,x+\phi_2) \right\}
\\&
=\D x{}\left\{\kappa({x+\phi_1},x+\phi_2)\D x{u_\phiv} \right\},
\end{align*}
namely the right-hand side of~\eqref{eq:shdifpde}.
Hence, provided \pde~\eqref{eq:emdifpde} has boundary conditions of \(\ell_i\)-peridicity in~\(x_i\), every solution of the embedding \pde~\eqref{eq:emdifpde} gives a solution of the original \pde~\eqref{eq:hdifpde} for every multi-dimensional phase-shift~\phiv\ of the heterogeneity.

In particular, the field~\(u_\ov(t,x)\) satisfies the given heterogeneous \pde~\eqref{eq:hdifpde}.
\end{proof}

\begin{lemma}[converse]\label{lemcon}
Suppose we have a set of solutions~\(u_\phiv(t,x)\) of the phase-shifted \pde~\eqref{eq:shdifpde}---a set parametrised by the phase vector~\phiv---and the set depends smoothly upon~\(t,x,\phiv\).
Then the field \(\fu(t,x,x_1,x_2):=u_{(x_1-x,x_2-x)}(t,x)\) satisfies the embedding \pde~\eqref{eq:emdifpde}.
\end{lemma}

\begin{proof} 
First, from the \pde~\eqref{eq:emdifpde}, consider
\begin{align*}
\D{x_0}\fu+\D {x_1}\fu+\D {x_2}\fu
&=-\D{\phi_1}{u_\phiv}-\D{\phi_2}{u_\phiv}+\D x{u_\phiv}+\D{\phi_1}{u_\phiv}+\D{\phi_2}{u_\phiv}
=\D x{u_\phiv}\,.
\end{align*}
Second, since \(\phi_i=x_i-x\)\,, that is \(x_i=x+\phi_i\)\,, then
\begin{align*}
\left(\D{x_0}{}+\D {x_1}{}+\D {x_2}{}\right)\kappa(x_1,x_2)
&=\D x{}\kappa(x+\phi_1,x+\phi_2).
\end{align*}
Thirdly, hence the right-hand-side of \pde~\eqref{eq:emdifpde} becomes
\begin{align*}
&\left(\D {x_0}{}+\D {x_1}{}+\D {x_2}{}\right)\left\{\kappa({x_1},x_2)\D x{u_\phiv}\right\}
\\&
=\D x{}\kappa(x+\phi_1,x+\phi_2)\D x{u_\phiv}
+\kappa(x_1,x_2)\left(\D {x_0}{}+\D {x_1}{}+\D {x_2}{}\right)\D x{u_\phiv}
\\&
=\D x{}\kappa(x+\phi_1,x+\phi_2)\D x{u_\phiv}
+\kappa(x+\phi_1,x+\phi_2)\DD x{u_\phiv}
\\&
=\D x{}\left\{\kappa(x+\phi_1,x+\phi_2)\D x{u_\phiv}\right\}
\,,
\end{align*}
the right-hand side of \pde~\eqref{eq:shdifpde}.
Lastly, since \(\D t\fu=\D t{u_\phiv}\) it follows that \(\fu:=u_{(x_1-x,x_2-x)}(t,x)\) satisfies the embedding \pde~\eqref{eq:emdifpde}.
\end{proof}

Corresponding results to \cref{lemeqv,lemcon} hold for the heterogeneous wave \pde\ \(u_{tt}=\partial_x\{\kappa(x)u_x\}\), with~\(\partial_t\) replaced by~\(\partial_{tt}\).

Consequently, \pde{}s~\eqref{eq:emdifpde} and~\eqref{eq:shdifpde} are equivalent, and they may provide us with a set of solutions for an ensemble of materials all with the same heterogeneity structure, but with the structural phase of the material shifted through all possibilities.
The key difference between \pde{}s~\cref{eq:emdifpde,eq:shdifpde} is that although~\eqref{eq:shdifpde} is heterogeneous in space~\(x\), the embedding \pde~\eqref{eq:emdifpde} is \emph{homogeneous} in~\(x_0\).
Because of this homogeneity, \cref{secsmh} is empowered to apply a rigorous theory for slow variations in space that leads to the \text{desired homogenized \pde.}

\section{The slow manifold of homogenization}
\label{secsmh}

Let's analyse the embedding \pde~\eqref{eq:emdifpde} for a useful slow manifold.  
This slow manifold expresses and supports the emergence of a precise homogenization of the original heterogeneous \pde~\eqref{eq:hdifpde}.
Since the \pde{}s herein are linear, the slow manifold is more specifically a slow subspace, but I use the term \emph{manifold} as the same framework and theory immediately also applies to cognate \text{nonlinear systems.}

Rigorous theory \cite[]{Roberts2013a} inspired by earlier formal arguments \cite[]{Roberts88a, Roberts96a} establishes how to create a \pde\ model for the macroscale spatial structure of \pde\ solutions in cylindrical domains~\cD\ like \eqref{fig:pattensemble}.
The technique is to base analysis on the case where variations in~\(x_0\) are approximately negligible, and then treat slow, macroscale, variations in~\(x_0\) as a regular perturbation \cite[Part~III]{Roberts2014a}.%
\footnote{Alternatively, in linear problems one could justify the analysis via a Fourier transform in~\(x_0\) \cite[\S2 and Ch.7 resp.]{Roberts88a, Roberts2014a}.  
However, for nonlinear problems it is better to analyse in physical space, so we do so herein---except for \cref{secsrh}.}

The analysis proceeds to any chosen order~\(N\) in the `small' derivatives in~\(x_0\).
For two examples, the usual leading order homogenizations are the case \(N=2\), and the so-called second-order homogenizations\footcite[e.g.,][]{Anthoine2010, Cornaggia2020, AewisHii2022} correspond to the higher-order \(N=4\).
To ensure the required derivatives and operators exist and have requisite properties, we analyse the systems in Sobolov spaces \(\HH_D:=W^{N+1,2}(D)\) for various spatial domains~\(D\).

To establish the basis of a slow manifold, consider \pde~\eqref{eq:emdifpde}, in the Sobolov space~\(\HH_\cD\), with \(\D{x_0}{}\) neglected:
\begin{equation}
\D t\fu=\left(\D {x_1}{}+\D {x_2}{}\right)\left\{\kappa(x_1,x_2)\left(\D {x_1}\fu+\D {x_2}\fu \right)\right\},
\label{eq:emdifpde0}
\end{equation}
to be solved with \(\ell_i\)-periodic boundary conditions in~\({x_i}\).
The basis applies at each and every~\(x_0\).
Hence at every~\(x_0\) we consider \pde~\eqref{eq:emdifpde0} on the Sobolov space~\(\HH_\cC\) where \(\cC\)~is the \(x_1x_2\)-cross-section of~\cD.

\paragraph{Equilibria}
A family of equilibria is \(\fu(t,x_1,x_2)=U\), constant in~\(x_1,x_2\), for every~\(U\). 
This family, in~\(\HH_\cC\) for every~\(x_0\), forms a subspace of equilibria~\EE.

\subsection{Spectrum at each equilibrium}
\label{secsee}

We explore the spectrum about every equilibria in~\EE\ as the spectrum is crucial in identifying the existence and properties of invariant manifolds.
Since \pde~\eqref{eq:emdifpde0} is linear, the perturbation problem is identical at every~\(U\), namely~\eqref{eq:emdifpde0} itself.

\begin{example}\label{egKappaConst}
Consider the case where the heterogeneity~\(\kappa\) is constant.  
In the vector space~\(\HH_\cC\), and upon defining wavenumbers \(k_i:=2\pi/\ell_i\)\,, a complete set of linearly independent eigenfunctions of the right-hand side of~\eqref{eq:emdifpde0} are \(\e^{\i(m k_1 x_1+n k_2 x_2)}\) for every integer~\(m,n\).
The spectrum of eigenvalues is then \(\lambda_{mn}:=-\kappa (mk_1+nk_2)^2\).
This spectrum has one zero-eigenvalue and the rest are negative: \(\lambda_{mn}\leq -\kappa k_1^2 <0\) (recall \(\ell_1>\ell_2\) so \(k_1<k_2\)).
\end{example}

For general heterogeneity~\(\kappa\), we here establish that the eigenvalue problem associated with~\eqref{eq:emdifpde0} is self-adjoint and hence has only real eigenvalues.
Denote the right-hand side operator as~\(\cL_0:=(\partial_{x_1}+\partial_{x_2})\big\{\kappa(x_1,x_2)(\partial_{x_1}+\partial_{x_2})\big\}\) to be solved with \(\ell_i\)-periodic boundary conditions in~\({x_i}\).
Also, we use the inner product \(\left<\fu,\fv\right>:=\iint_{\cC} \fu\bar\fv \d{x_1}\d{x_2}\) over the rectangular cross-section \(\cC:=[0,\ell_1]\times[0,\ell_2]\), where an overbar denotes complex conjugation.
Then, for every pair of functions~\(\fu,\fv\in\HH\), consider
\begin{align*}
\left<\cL_0\fu,\fv\right>&
=\iint_{\cC} \bar\fv(\partial_{x_1}+\partial_{x_2})\big\{\kappa (\fu_{x_1}+\fu_{x_2})\big\} \d{x_1}\d{x_2}
\\&
=\iint_{\cC} \bar\fv\partial_{x_1}\big\{\kappa (\fu_{x_1}+\fu_{x_2})\big\} + \bar\fv\partial_{x_2}\big\{\kappa (\fu_{x_1}+\fu_{x_2})\big\} \d{x_1}\d{x_2}
\\&
=\int_0^{\ell_2} \big[\bar\fv\kappa (\fu_{x_1}+\fu_{x_2})\big]_0^{\ell_1}
-\int_0^{\ell_1} \bar\fv_{x_1}\kappa (\fu_{x_1}+\fu_{x_2})  \d{x_1}\d{x_2}
\\&\quad{}
+ \int_0^{\ell_1} \big[\bar\fv\kappa (\fu_{x_1}+\fu_{x_2})\big]_0^{\ell_2}
 -\int_0^{\ell_2} \bar\fv_{x_2}\kappa (\fu_{x_1}+\fu_{x_2}) \d{x_2}\d{x_1}
\\&
=-\iint_{\cC} (\bar\fv_{x_1}+\bar\fv_{x_2})\kappa (\fu_{x_1}+\fu_{x_2}) \d{x_1}\d{x_2}
\end{align*}
as the boundary contributions vanish due to the \(\ell_i\)-periodicity.
Reverse the above steps with the roles of~\fu\ and~\fv\ swapped, to then deduce \(\left<\cL_0\fu,\fv\right> =\left<\fu,\cL_0\fv\right>\).
Thus the operator~\(\cL_0\)\ is self-adjoint.

Consequently every eigenvalue is real, and the eigenvectors are pairwise orthogonal (or may be chosen to be so).

The next step is to establish that all the non-zero eigenvalues of~\(\cL_0\)\ are negative and bounded away from zero.
Consider the eigen-problem \(\cL_0\fv=\lambda\fv\) for cross-sectional structures \(\fv\in \HH_\cC\).
Apply \(\left<\cdot,\fv\right>\) to the eigen-problem to see
\begin{align*}
\lambda\left<\fv,\fv\right>
&=\iint_\cC \bar\fv(\partial_{x_1}+\partial_{x_2})\big\{\kappa (\fv_{x_1}+\fv_{x_2})\big\}\d{x_1}\d{x_2}
\\&
=\cdots=-\iint_\cC \kappa |\fv_{x_1}+\fv_{x_2}|^2\d{x_1}\d{x_2}
\\&
\leq-\kappa_{\min}\iint_\cC  |\fv_{x_1}+\fv_{x_2}|^2\d{x_1}\d{x_2}
\,,
\end{align*}
via integration by parts, and the \(\ell_i\)-periodic conditions in~\(x_i\).
To deduce an upper bound on the non-zero eigenvalues~\(\lambda\), let's minimise \(\iint_\cC  |\fv_{x_1}+\fv_{x_2}|^2\d{x_1}\d{x_2}\) subject to \(\left<\fv,\fv\right>=\iint_\cC|\fv|^2\d{x_1}\d{x_2}=1\) and that \(\fv({x_1},x_2)\)~has zero mean (being orthogonal to the connstant eigenvector corresponding to \(\lambda=0\)).
Introduce a Lagrange multiplier~\(\mu\), and use calculus of variations to minimise \(\iint_\cC  |\fv_{x_1}+\fv_{x_2}|^2 \d{x_1}\d{x_2} +\mu\left( \iint_\cC|\fv|^2 \d{x_1}\d{x_2}-1 \right)\).
Via integration by parts, it follows that a minimiser must satisfy \((\partial_{x_1}+\partial_{x_2})^2\fv=-\mu\fv\) over \(\fv\in\HH_\cC\).
This is the problem \cref{egKappaConst} solved, giving potential minimisers as~\(\e^{\i(m k_1 x_1+n k_2 x_2)}\) for multiplier \(\mu=(mk_1+nk_2)^2\).
The zero mean condition requires at least one of~\(m,n\) non-zero, and since \(\ell_1>\ell_2\), so minima are thus linear combinations of~\(\e^{\i m k_1 x_1}\) for \(m=\pm1\).
These give \(\iint_\cC  |\fv_{x_1}+\fv_{x_2}|^2\d{x_1}\d{x_2}/\left<\fv,\fv\right>=k_1^2=4\pi^2/\ell_1^2\).
Hence we obtain the upper bound on the non-zero eigenvalues of \(\lambda\leq-\beta_1\) for \(\beta_1:=4\pi^2\kappa_{\min}/\ell_1^2\).

In the case of the heterogeneous wave propagation problem, \(u_{tt}=\partial_x\{\kappa(x)u_x\}\), this subsection's eigenvalue analysis establishes that for the cross-sectional dynamics there exists a zero eigenvalue separated from pure imaginary eigenvalues of magnitude \text{at least \(\sqrt{\beta_1}=2\pi\sqrt{\kappa_{\min}}/\ell_1\).}

\subsection{Results of slowly varying theory}
\label{secrcvt}

Let's comment the preconditions for the next \cref{prosvpde} as listed in  Assumption~2 of \cite{Roberts2013a}. 
First, \cref{secsee} establishes that the spectrum of~\(\cL_0\) is as required for the Hilbert space~\(\HH_\cC\).
Hence we decompose the space into two closed \(\cL_0\)-invariant subspaces, the centre~\(\EE_c^0\) and the stable~\(\EE_s^0\).
Since the operator~\(\cL_0\) is self-adjoint it generates the required strongly continuous semigroups in the Sobolev space~\(\HH_\cC\) \cite[e.g.,][Ch.6]{Carr81}.
Hence the following result holds for homogenization.

\begin{proposition}[homogenized \pde]\label{prosvpde}
By Proposition~6 of \cite{Roberts2013a}, for every chosen truncation order~\(N\), in a domain~\XX\ of `slowly varying solutions', and after the decay of transients~\Ord{\e^{-\beta't}} for \(\beta'\approx\beta_1\), the mean field \(U(t,x):=(\ell_1\ell_2)^{-1} \iint_\cC \fu(t,x,x_1,x_2) \d{x_1}\d{x_2}\) satisfies a generalised homogeneous \pde
\begin{equation}
\D tU=\sum_{n=0}^N K_n\Dn xnU\,, \quad x\in\XX\subset[0,L],
\label{eq:genhomopde}
\end{equation}
in terms of some deterministic constants~\(K_n\), to a quantifiable error \cite[(9)--(10) and~(23), resp.]{Roberts2013a}.
\end{proposition}
The above mentioned error~(23) \cite[from][]{Roberts2013a} involves too many new parameters to meaningfully detail here---but typically the most important factor at each and every location \(x=X\) is the \((N+1)\)th~derivative of~\(\fu\) in a spatial neighbourhood of~\(X\).
The open subset~\XX\ of the domain~\([0,L]\) is then that part of the domain where the error in~\eqref{eq:genhomopde} \cite[eqn.~(23)]{Roberts2013a}, dominantly~\(\partial_x^{N+1}\fu\), is small enough for the modelling purposes at hand.
Also, using the specific mean field defined in \cref{prosvpde}---the usual `cell average'---is optional: one may instead parametrise the macroscale slow manifold in terms of any reasonable alternative amplitude \cite[e.g.,][\S5.3.3]{Roberts2014a}.

A generating function argument \cite[Corollaries 12 and~13]{Roberts2013a} establishes that classic, formal, `slowly varying' analyses of the embedding \pde~\eqref{eq:emdifpde} are valid.
One just needs to formally treat derivatives in the `large' direction,~\(\D{x_0}{}\), as `small' in a consistent asymptotic sense {\cite[e.g.,][]{Roberts88a, Roberts96a}}.

In the case of the heterogeneous wave propagation problem, \(u_{tt}=\partial_x\{\kappa(x)u_x\}\), I conjecture that a backward theory \cite[]{Roberts2018a, Hochs2019} approach to the slowly varying theory here would support a cognate version of \cref{prosvpde} about the generalised homogenous wave \pde\ \(U_{tt}=\sum_{n=0}^N K_n\partial_x^nU\) describing its dispersive macroscale waves.

\section{Example construction of homogenized PDE} 
\label{secxchm}

Corollary~13 of \cite{Roberts2013a} established that the procedure of \cite{Roberts88a}, previously viewed as formal, is indeed a rigorous way of constructing the slow manifold modelling \pde{}s~\eqref{eq:genhomopde}.
Here we use the computer algebra system (\textsc{cas}) version of the procedure \cite[]{Roberts96a}, as detailed in generality and examples of a book \cite[Part~III]{Roberts2014a}.

\cref{seccachm} lists and documents the computer algebra code.
The code constructs the homogenization for the case of heterogeneity 
\begin{equation}
\kappa(x_1,x_2):=1/[1+a_1\cos k_1x_1+a_2\cos k_2x_2)]
\label{egkappa}
\end{equation}  
in terms of microscale wavenumbers \(k_i:=2\pi/\ell_i\)
for the given microscale periodicities~\(\ell_1,\ell_2\).
We write approximations to the slow manifold model of the embedding \pde~\eqref{eq:emdifpde} in terms of a `mean' field~\(U(t,x)\) that evolves according to an homogenized \pde\ of the form~\eqref{eq:genhomopde}.

\paragraph{Quickly verify an approximation}
First, \cref{seccaqvg} verifies that an approximate field is 
\begin{subequations}\label{eqsqvg}%
\begin{align}
\fu(t,x,x_1,x_2)&=U
+(a_1/k_1\sin k_1x_1 +a_2/k_2\sin k_2x_2)\D xU
\nonumber\\&\quad{}
+(a_1/k_1^2\cos k_1x_1 +a_2/k_2^2\cos k_2x_2)\DD xU+\cdots
\label{eqqvgu}
\end{align}
such that the mean field~\(U\) evolves according to the homogenized \pde\  
\begin{equation}
\D tU=\DD xU+\cdots
\label{eqqvgg}
\end{equation}
\end{subequations}
(the coefficient~\(1\) of the mean diffusivity is the classic harmonic mean of the specific heterogeneous diffusivity~\eqref{egkappa}).
To verify, substitute~\eqref{eqsqvg} into the governing \pde~\eqref{eq:emdifpde} and find the \pde's residual contains only terms in~\(U_{xxx}\) and higher derivatives.
We denote this by saying that the residual is zero to an error~\Ord{\partial_x^3}.
Since the \pde\ residual is~\Ord{\partial_x^3}, then theory \cite[Prop.~6]{Roberts2013a} assures us that the slow manifold~\eqref{eqsqvg} is correct to error~\Ord{\partial_x^3}.

\subsection{Iteration systematically constructs}
\label{secisc}

To construct approximations systematically, we repeatedly compute the residual, which then drives corrections, until the residual is zero to the specified order of error.  
Theory then assures us the slow manifold is approximated to the same order of error \cite[]{Roberts2013a}.
For example, the \cas\ iteratively constructs the improved homogenization (the case \(N=6\)) that
\begin{align}
\D tU&=\DD xU
+\tfrac12(a_1^2/k_1^2+a_2^2/k_2^2)\Dn x4U
\nonumber\\&\quad{}
+\left[ \tfrac12(a_1^2/k_1^2+a_2^2/k_2^2)^2 
        -2(a_1^2/k_1^4+a_2^2/k_2^4) \right]\Dn x6U
+\Ord{\partial_x^7}.
\label{eqimphpde}
\end{align}
The code of \cref{seccaisc} constructs this generalised homogenization~\eqref{eqimphpde}, and executes in about 26~seconds.
In the case of the heterogeneous wave propagation problem, \(u_{tt}=\partial_x\{\kappa(x)u_x\}\), we expect the macroscale homogenized dispersive wave \pde\ to have \(U_{tt}\)~equal to the same right-hand side as the above.

In use, any specific truncation of higher-order homogenization, such as~\eqref{eqimphpde}, may need some asymptotically consistent regularisation \cite[e.g.,][\S2]{Benjamin1972}.  
For an example of regularisation, let's suppose we choose the \(N=4\) truncation \(U_t=U_{xx}+\alpha^2U_{xxxx}\) with parameter \(\alpha:= \sqrt{(a_1^2/k_1^2+a_2^2/k_2^2)/2}\).
The fourth derivative term~\(+\alpha^2U_{xxxx}\) is undesirably destabilising.
But to the same order of error, namely~\Ord{\partial_x^6}, the chosen truncation is asymptotically equivalent to the regularised equation \(U_t=(1-\alpha^2\partial_x^2)^{-1} U_{xx}\).
This regularisation may be alternatively written as
\begin{equation}
(1-\alpha^2\partial_x^2)U_t=U_{xx}\,,
\quad\text{or}\quad
U_t =\tfrac1{2\alpha} \e^{-|x|/\alpha}\star U_{xx}\,.
\label{eqregularised}
\end{equation}
The spatially nonlocal convolution~\(\e^{-|x/\alpha|}\star\), whether explicit or implicit, stabilises the regularised chosen truncation while preserving asymptotic consistency \cite[cf.,][\S2.2]{Mercer94a}.  
Such nonlocal \pde{}s have previously been found desirable in various  homogenizations of heterogeneous systems \cite[e.g.,][pp.2--3]{Cooper2018}.

\paragraph{Optional} We may set \(a_2=0\) to get simpler results for a single periodic component in the heterogeneity, namely \(\kappa(x_1)=1/[1+a_1\cos k_1x_1]\).
In this case it becomes feasible to construct the homogenization to high-order, \(N=34\)\,, in spatial derivatives.
To construct expressions that are exact in heterogeneity amplitude~\(a_1\), we also need to construct to the same order in~\(a_1\) because the evidence is that here the expansions in homogeneity amplitude~\(a_1\) truncate at just less than the same order as the order of derivatives.
For this option (\cref{seccaisc}), the algebra of order \(N=34\) takes roughly one minute compute time.

\paragraph{In general}
The code of \cref{seccaisc} assumes the heterogeneity \(\kappa\approx 1\), specifically \(\kappa=1/[1+\kappa']\) for some~\(\kappa'\propto a\).  
Then multiplication by~\(\kappa\) is realised as multiplication by \(\sum_{n=0}^{N-1}(-\kappa')^n\).

The method \cite[]{Roberts96a} is to start the iteration with the leading `trivial' approximation for the field and its evolution: that is,  \(u\approx u^{(1)}:= U\) such that \(\D tU\approx g^{(1)}:=0\).
The \(n\)th~iteration is given the approximations~\(u^{(n)}\) and~\(g^{(n)}\), and starts by evaluating, via the flux \(f^{(n)}:=-\kappa(u^{(n)}_x+u^{(n)}_{x_1}+u^{(n)}_{x_2})\), the residual of the embedding \pde~\eqref{eq:emdifpde}:
\begin{equation*}
\Res^{(n)}:=u^{(n)}_t+f^{(n)}_x+f^{(n)}_{x_1}+f^{(n)}_{x_2}\,.
\end{equation*}
First update the evolution via the solvability integral that here is the mean over~\(x_1,x_2\):  \(g^{(n+1)}:=g^{(n)}+g'\) where \(g':=\overline{\Res^{(n)}}\).
Second, update the field by a correction~\(u'\) through using two steps to solve \((\partial_{x_1} +\partial_{x_2}) \big[ \kappa (u'_{x_1} +u'_{x_2}) \big] =\Res^{(n)} +g'\): 
step~1 integrates to solve \(v_{x_1}+v_{x_2}=\Res^{(n)} +g'\) with zero mean to avoid unbounded updates;
and step~2 integrates again to solve \(u'_{x_1}+u'_{x_2}=v/\kappa\) with zero mean.
This second ``zero mean'' ensures preservation of our \emph{choice} that the macroscale field parameter~\(U\) is to be the local mean of the microscale field~\(u\).
Then update the field approximation with \(u^{(n+1)}:=u^{(n)}+u'\).

The iterative loop terminates when the \pde\ residual~\(\Res^{(n)}\) is zero to the specified orders of error.
Then the constructed homogenization, such as~\cref{eqsqvg,eqimphpde}, is correct to the same order \cite[Prop.~6]{Roberts2013a}.

\paragraph{Analogy with machine learning}  
The recursive iteration used in the above construction was originally developed nearly 30~years ago \cite[]{Roberts96a}.  
Let's draw an analogy between this iteration and machine learning algorithms where an AI learns the generic form of the macroscale evolution from many thousands of simulations\footcite[e.g.,][]{GonzalezGarcia98, Frank2020,  Linot2020}.
In the algorithm here, each recursive iteration is analogous to a layer in a deep neural network, evaluating residuals is analogous to a nonlinear neuronal function, the linear update corrections are analogous to using weighted linear combination of outputs of one layer as the inputs of the next layer.  
But here the recursion is a `smart neural network' in that the neurones and the `linear weights' are crafted to the problem at hand using the physics encoded in the \pde\ operator.
Further, being analytic, a single analysis encompasses all `data points' in the state space's domain (here all `slowly varying' functions in~\(\HH_\cD\)), not just a finite sample as in machine learning (possibly a biased sample).
Consequently, via such algorithms, I contend that mathematicians have for decades been using such recursive iteration for implementing smart analogues of machine learning.
Such analytic learning empowers the physical interpretation, validation, and verification required by modern science \cite[e.g.][]{Brenner2021}.

\section{Spatial resolution of homogenization}
\label{secsrh}

An homogenized \pde, such as~\eqref{eq:genhomopde,eqimphpde}, accurately predicts the evolution of the mean field~\(U(t,x)\) \emph{provided} that, on the scale of the heterogeneity, the spatial gradients are sufficiently small.  
That is, provided the solutions are sufficiently slowly varying in space.
\cite{Matous2017} [\S3.2] phrased the proviso as ``the scale of the microstructural fluctuations,~[\(\ell_1\)], must be \ldots\ much smaller than the macroscopic field fluctuations,~\(\ell_M\)'' \cite[expressed similarly by][\S2.2.2]{VinhPhuNguyen2011}.
Question: what does ``sufficiently small'', ``sufficiently slowly varying'', or ``much smaller'' actually mean physically?
This section provides an indicative \text{\emph{quantitative} answer.}

Here we show that in a prototypical family of homogenization problems, we only need to require that the 
\begin{equation}
\text{macroscale lengths } \ell_M\gtrsim2\ell_1\,.
\label{eqMacro}
\end{equation}
This lower bound on the macroscale lengths is approximate and indicative because the actual performance of an homogenization depends upon the details of the heterogeneity~\(\kappa\), a chosen truncation of the homogenized \pde, any specific regularisation of the \pde, the desired error of predictions, and so on. 
But the numerical coefficient~\(2\) in~\eqref{eqMacro} arises from a well-established \text{quantitative procedure.}

The procedure we invoke is that of high-order construction followed by Domb--Sykes plots \cite[e.g.,][]{Domb57, Vandyke84} and Mercer--Roberts plots \cite[Appendix]{Mercer90} that estimate the distance to the nearest, convergence limiting, singularity in series approximations.
We analyse an homogenized \pde~\eqref{eq:genhomopde,eqimphpde} by considering its corresponding Fourier transform (\(\partial_x\leftrightarrow\i k\)): \(\D t{\tilde U}=\sum_{n=0}^\infty (\i k)^nK_n\tilde U\) in terms of the Fourier transform~\(\tilde U(t,k)\) of the mean field~\(U(t,x)\).  
Then the coefficient \(\cK(k):=\sum_{n=0}^\infty (\i^nK_n)k^n\) is a series in wavenumber~\(k\).
If we can show convergence for \(|k|<k_*\) for some wavenumber bound~\(k_*\), then we deduce that the homogenized \pde\ is valid for spatial structures of length scale \(\ell_M>2\pi/k_*\).
We know that the radius of convergence~\(k_*\) is the distance to the singularity of~\(\cK(k)\) nearest the origin (\(k=0\)).  
But we only know~\(N\) terms in the series for~\(\cK(k)\), so Domb--Sykes and Mercer--Roberts plots estimate the distance~\(k_*\) by extrapolating the `ratio test' from finite~\(n\) to infinity, in the cases of a single singularity or a complex conjugate pair of singularities, respectively.
For example, this technique has been used to quantitatively predict the spatial limit on \pde{}s modelling shear dispersion in channels and pipes \cite[]{Mercer90, Mercer94a}.
\begin{figure}
\centering
\caption{\label{fig:homoRadius}%
macroscale wavenumbers~\(k\) of convergence are estimated from Mercer--Roberts plots \protect\cite[Appendix]{Mercer90} for the series~\(\cK(k)\) in the case of heterogeneity~\eqref{eqDiff1}, here with \(a=0.975\).  
We use the first~17 terms of the series in~\(k^2\).  
Extrapolation to \(1/n=0\) estimates the location of a pair of complex conjugate singularities that limit the radius \text{of convergence.}}
%
\tikzsetnextfilename{Figs/homoRadius}
\definecolor{mycolor1}{rgb}{0.13359,0.48867,0.73828}%
\definecolor{mycolor2}{rgb}{0.46758,0.53789,0.00000}%
\begin{tikzpicture}

\begin{axis}[name=top,%
xmin=0,
xmax=0.25,
xlabel style={font=\color{white!15!black}},
xlabel={$1/n$},
ymin=2.8,
ymax=3.8,
ylabel style={font=\color{white!15!black}},
ylabel={$1/k_*^2$},
axis background/.style={fill=white},
width=100mm,height=40mm
]
\addplot [color=mycolor1, dotted, mark=o, mark options={solid, mycolor1}, forget plot]
  table[row sep=crcr]{%
0.25	3.50139858386572\\
0.2	3.26518508995875\\
0.166666666666667	3.16899588685437\\
0.142857142857143	3.11855157366568\\
0.125	3.11116556288059\\
0.111111111111111	3.15343459085316\\
0.1	3.23515130526846\\
0.0909090909090909	3.32192259938042\\
0.0833333333333333	3.37862463490489\\
0.0769230769230769	3.39776559499579\\
0.0714285714285714	3.39722033691974\\
0.0666666666666667	3.40004684742514\\
};
\addplot [color=mycolor2, forget plot]
  table[row sep=crcr]{%
0.25	2.8154840202015\\
0.2	2.97602265543671\\
0.166666666666667	3.08304841226018\\
0.142857142857143	3.15949538141981\\
0.125	3.21683060828952\\
0.111111111111111	3.26142467363264\\
0.1	3.29709992590713\\
0.0909090909090909	3.32628876867717\\
0.0833333333333333	3.35061280431886\\
0.0769230769230769	3.37119468063107\\
0.0714285714285714	3.38883628889868\\
0.0666666666666667	3.4041256827306\\
0	3.61817719637755\\
};
\end{axis}

\begin{axis}[at={(top.below south west)},anchor={north west},%
xmin=0,
xmax=0.08,
xlabel style={font=\color{white!15!black}},
xlabel={$1/n^2$},
ymin=0.88,
ymax=1,
ylabel style={font=\color{white!15!black}},
ylabel={$\cos2\theta_*$},
axis background/.style={fill=white},
width=100mm,height=40mm
]
\addplot [color=mycolor1, dotted, mark=o, mark options={solid, mycolor1}, forget plot]
  table[row sep=crcr]{%
0.0625	0.978229282637076\\
0.04	0.959866923491973\\
0.0277777777777778	0.948073401277267\\
0.0204081632653061	0.935152377552595\\
0.015625	0.920120570704884\\
0.0123456790123457	0.905296380937195\\
0.01	0.895237565265754\\
0.00826446280991736	0.892583676936205\\
0.00694444444444444	0.894856958394094\\
0.00591715976331361	0.897318577105806\\
0.00510204081632653	0.897112987627055\\
0.00444444444444444	0.894200513122634\\
};
\addplot [color=mycolor2, forget plot]
  table[row sep=crcr]{%
0.0625	0.997260495208701\\
0.04	0.956956751197994\\
0.0277777777777778	0.935063359389709\\
0.0204081632653061	0.921862334644453\\
0.015625	0.913294361853062\\
0.0123456790123457	0.907420187914601\\
0.01	0.903218425850385\\
0.00826446280991736	0.900109597111267\\
0.00694444444444444	0.897745077898313\\
0.00591715976331361	0.895904925951006\\
0.00510204081632653	0.894444821712\\
0.00444444444444444	0.893266884119346\\
0	0.885305650734515\\
};
\end{axis}
\end{tikzpicture}%
\end{figure}
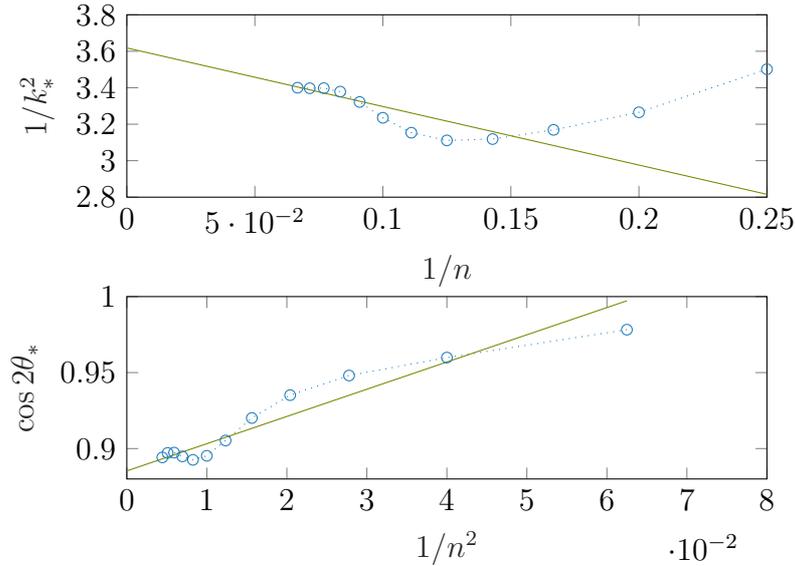%
\cref{fig:homoRadius} is an example Mercer--Roberts plot for one homogenization case developed here.  
The extrapolations to \(1/n=0\) (\(n\to\infty\)) in the figure indicate that the convergence is limited by a complex conjugate pair of singularities at a distance in wavenumber space \text{of \(k_*\approx 1/\sqrt{3.6} \approx 0.52\).}

To obtain such plots and estimates we need to compute the homogenization series~\eqref{eq:genhomopde,eqimphpde} to high-order.
For the case of quasi-periodic heterogeneity the algebraic complexity explodes combinatorially with increasing order.  
Consequently, as a representative example, we explore the accessible case of homogenization with a single microscale period, specifically the heterogeneity \begin{equation}
\kappa(x):=1/[1+a\cos x].
\label{eqDiff1}
\end{equation}
The computer algebra code of \eqref{secxchm} constructs the homogenization in this case via the option of setting \(a_2=0\), \(a_1=a\), and \(k_1=1\).
In this case the analysis is non-dimensionalised on the microscale length so that non-dimensionally \(\ell_1=2\pi\).
Executing the code with this option gives that the homogenized \pde~\eqref{eq:genhomopde} begins with
\begin{align}
U_t&=\big[\partial_x^2 
+\tfrac12a^2\partial_x^4 
+(-2a^2-\tfrac12a^4)\partial_x^6
+(8a^2-\tfrac{135}{32}a^4+\tfrac58a^6)\partial_x^8
\nonumber\\&\qquad{}
+(-32a^2+\tfrac{879}{32}a^4-\tfrac{261}{32}a^6+\tfrac78a^8) \partial_x^{10}
+\Ord{\partial_x^{12}}\big]U.
\label{eqHoHomo1}
\end{align}
As a first step in investigating such high-order homogenization, observe that the \(a^2\)~terms are evidently all of the form \(\tfrac12a^2 \partial_x^4 (-4\partial_x^2)^{n/2-2}\) (this form holds up to at least order \(n=34\)).
Hence in Fourier space (\(\partial_x\leftrightarrow\i k\)), the \pde~\eqref{eqHoHomo1} becomes the series for \(\tilde U_t=\big[-k^2+\tfrac12a^2k^4/(1-4k^2)+\Ord{a^4}\big]\tilde U\).
The divisor~\((1-4k^2)\) in this expression indicates that, for small heterogeneity amplitude~\(a\), the wavenumbers resolvable by such an homogenization series are limited by the pole singularities at \(k=\pm \tfrac12\), that is, a bound on the convergence of the series is \(k_*=1/2\).
Hence, for small~\(a\), the homogenized \pde~\eqref{eqHoHomo1} potentially resolves all wavelengths longer than \(\ell_*:=2\pi/k_*=4\pi=2\ell_1\).%
\footnote{The limit \(a\to0\) is non-uniform.  
At \(a=0\) there is no heterogeneity and the `homogenized' \pde~\eqref{eq:edifpde} is precisely the original \pde~\eqref{eq:hdifpde} and therefore valid for \emph{all} wavelengths and wavenumbers.  
Whereas for small non-zero~\(a\) the singularities at wavenumber \(k=\pm\tfrac12\) persist and limit validity for small~\(a\)---it is just that the singularities have a strength~\Ord{a^2} that vanishes \text{as \(a\to0\).}}

I conjecture that the singularities at wavenumber~\(k=\pm\tfrac12\) are connected to, in wave propagation through heterogeneous media, the 
well known phenomenon of spectral gaps, namely frequencies at which no wave can propagate through the underlying medium \cite[e.g.,][p.2]{Cooper2018}.

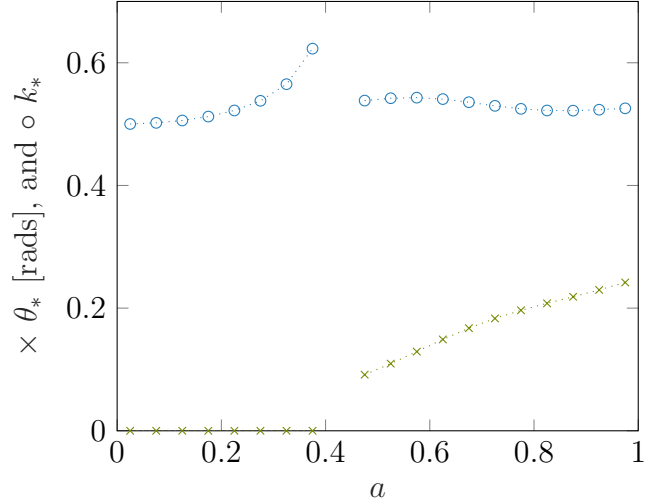
\begin{figure}
\centering
\begin{tabular}{@{}cc@{}}
\parbox[t]{45mm}{\caption{\label{fig:homoRadii}%
In complex wavenumber space, plot estimates of the radius~\(k_*\) (circles) and angle~\(\theta_*\) (crosses) of the convergence limiting singularity of the homogenization~\eqref{eqHoHomo1} as a function of the heterogeneity amplitude~\(a\) in~\eqref{eqDiff1}.
The estimates around \(a\approx 0.4\) are unreliable (see text).
}}
&
\raisebox{-\height}{
%
\tikzsetnextfilename{Figs/homoRadii}
\definecolor{mycolor1}{rgb}{0.13359,0.48867,0.73828}%
\definecolor{mycolor2}{rgb}{0.46758,0.53789,0.00000}%
\begin{tikzpicture}

\begin{axis}[%
unbounded coords=jump,
xmin=0,
xmax=1,
xlabel style={font=\color{white!15!black}},
xlabel={$a$},
ymin=0,
ymax=0.7,
ylabel style={font=\color{white!15!black}},
ylabel={$\times\ \theta_*$ [rads],   and   $\circ\ k_{*}$},
axis background/.style={fill=white}
]
\addplot [color=mycolor1, dotted, mark=o, mark options={solid, mycolor1}, forget plot]
  table[row sep=crcr]{%
0.025	0.50022368863501\\
0.075	0.502049386419135\\
0.125	0.505908620113722\\
0.175	0.51230031497543\\
0.225	0.522253677398413\\
0.275	0.537991474841619\\
0.325	0.565124317310859\\
0.375	0.62304863154647\\
0.425	nan\\
0.475	0.538381016848204\\
0.525	0.542194396596626\\
0.575	0.54317647313786\\
0.625	0.540796958041994\\
0.675	0.535742324279963\\
0.725	0.529796527334537\\
0.775	0.524899318990965\\
0.825	0.522235670195735\\
0.875	0.521970428503884\\
0.925	0.523475830153696\\
0.975	0.525720707784815\\
};
\addplot [color=mycolor2, dotted, mark=x, mark options={solid, mycolor2}, forget plot]
  table[row sep=crcr]{%
0.025	0\\
0.075	0\\
0.125	0\\
0.175	0\\
0.225	0\\
0.275	0\\
0.325	0\\
0.375	0\\
0.425	nan\\
0.475	0.0914555312685378\\
0.525	0.109524995772015\\
0.575	0.1291102843037\\
0.625	0.148903694592694\\
0.675	0.167308088191274\\
0.725	0.183188123700781\\
0.775	0.19641539955337\\
0.825	0.207783631743032\\
0.875	0.218499281413115\\
0.925	0.229647306051174\\
0.975	0.241822713291633\\
};
\end{axis}
\end{tikzpicture}%
}%
\end{tabular}
\end{figure}%
\cref{fig:homoRadii}, for small~\(a\), confirms this bound on the wavenumbers and wavelengths.
The estimates plotted in \cref{fig:homoRadii} are obtained from \cref{secxchm} constructing the series~\eqref{eqHoHomo1} to errors of order~\(35\) in~\(\partial_x\), or equivalently to errors~\Ord{k^{35}}.
This \cas\ construction provides us with the first~17 non-zero coefficients in the series in~\(k^2\) for~\(\cK(k)\).
\begin{itemize}
\item For heterogeneity amplitude~\(a\lesssim0.4\) these coefficients have the same sign.
Consequently, Domb--Sykes plots \cite[e.g.,][]{Domb57, Vandyke84}, akin to the top panel of \cref{fig:homoRadius},  predict convergence limiting singularities at real \(k=\pm k_*\) as shown in the left-part of \cref{fig:homoRadii}.

\item For heterogeneity amplitude~\(a\gtrsim0.4\) the coefficients have a more complicated sign pattern, indicating the convergence limiting singularity is now a complex conjugate pair in the complex \(k^2\)-plane.
Mercer--Roberts plots \cite[Appendix]{Mercer90}, such as \cref{fig:homoRadius}, estimate the radius and angle of these convergence limiting singularities as shown in the right-part \text{of \cref{fig:homoRadii}.}
\end{itemize}

The two distinct behaviours summarised in \cref{fig:homoRadii} are likely due to singularities moving as a function of heterogeneity parameter~\(a\).  
We conjecture that for \(a\lesssim0.3\) there is a second real singularity for wavenumbers larger than the one at~\(k_*\).
As \(a\)~increases, both these singularities move and appear to collide for \(a\approx 0.35\).
After collision they separate and move out into the complex~\(k\) plane (\(a\gtrsim0.4\)).
The Domb--Sykes and Mercer--Roberts estimates are both unreliable near such collision, so the estimates for \(a\approx 0.3\) to~\(\approx0.4\) are unreliable.
With this caveat, the evidence of \cref{fig:homoRadii} indicates that across all~\(0<a<1\), our homogenization can resolve all wavenumbers \(|k|<1/2\).  
That is, high-order homogenization can resolve wavelengths longer than just twice the wavelength of the underlying periodicity. 

In the case of the heterogeneous wave propagation problem, \(u_{tt}=\partial_x\{\kappa(x)u_x\}\), we expect the same restrictions on the spatial resolution because we expect the right-hand side of \pde~\eqref{eqimphpde} to remain the same for the homogenized waves.

\paragraph{The case of quasi-periodic heterogeneity}
In the case of the quasi-periodic heterogeneity~\eqref{egkappa}, observe in the homogenization \pde~\eqref{eqimphpde} that the coefficients in~\(a_i\) are divided by wavenumbers~\(k_i\), and with increasing powers of~\(k_i\) as the order increases.
This indicates that higher-order terms in the series are dominated by the smallest wavenumber~\(k_i\), namely~\(k_1\).
That is, the higher-order terms should be dominated by the longest microscale heterogeneity, here that of \(a_1\cos(k_1x)\).
Hence we expect the limiting wavenumber and wavelengths to be the same as that for the single periodic heterogeneity~\eqref{eqDiff1}.
Consequently,  expect the bound~\eqref{eqMacro} to be valid for this quasi-periodic heterogeneity.

\paragraph{Temporal resolution}
The above exploration finds that the spatial resolution of the homogenization appears to vary little with heterogeneity strength~\(a\).
The temporal resolution is quite a different issue, and it may depend significantly upon~\(a\).
For example, \cref{fig:homoRadius} comes from the case \(a=0.975\) where the microscale heterogeneity varies by a factor of~\(80\) across a microscale period.  
In similar scenarios, the time-scale of decay to the homogenized model, roughly~\(1/\beta_1=\ell_1^2/(4\pi^2\kappa_{\min})\) (end of \cref{secsee}), may be several orders of magnitude longer than for small heterogeneity.
For valid predictions by the homogenization, the macroscale time scale of interest needs to be longer than that of this decay.

In the case of the heterogeneous wave propagation problem, \(u_{tt}=\partial_x\{\kappa(x)u_x\}\), the slow manifold is no longer attractive (unless there is some perturbing dissipative mechanism not included in this wave \pde).  
Instead we expect that the slow manifold acts as a `centre' for any fast oscillations arising from the initial conditions.
Indeed, conjectured backward theory \cite[e.g.,][\S2.5]{Roberts2018a} would more precisely assert that there is a system close to the \pde\ \(u_{tt}=\partial_x\{\kappa(x)u_x\}\) (close to some specified order in some sense) that possesses a homogeneous slow manifold for the macroscale that acts as a centre for all nearby dynamics.
In this way the backward theory would establish accurate modelling over ``timescales~\Ord{\varepsilon^{-\alpha}}'' \cite[p.3]{Abdulle2020c}, typically for exponent \(\alpha=N\), the order of truncation.

\section{Three-scale homogenization}
\label{sec3sh}

There is interest in physical systems possessing a spatial microscale, a mesoscale, and a macroscale\footcite[e.g.,][]{RamirezTorres2018, CruzGonzalez2022}.
Such three-scale systems are encompassed in the framework of \cref{secpse} when \(\ell_2\ll\ell_1\ll L\)\,. 
That is, the cylindrical domain of \cref{fig:pattensemble} becomes like a physical ruler in that the domain is significantly thinner vertically than its width, and its width is significantly thinner than its length.
Although \cref{secsrh} finds that these length ratios need to be at least two, they are not necessarily much bigger. 

In this scenario, the arguments of \cref{secpse} still apply: we solve the phase-shift embedding \pde~\eqref{eq:emdifpde} to obtain solutions of the heterogeneous \pde{}s~\eqref{eq:shdifpde,eq:hdifpde}.
It is the slow manifold support and analysis that differs in detail.

Because this rigorous homogenization could in principal be exact at finite scale separation, we find here the expected result that a slow manifold homogenization of a `meso-slow' manifold homogenization of the microscale system is the same as the slow manifold homogenization of the microscale (\cref{secsmh,secxchm}).
That is, this approach to homogenization is transitive \text{in principal.}

The same transitivity of homogenization would also hold in the case of the heterogeneous wave propagation problem, \(u_{tt}=\partial_x\{\kappa(x)u_x\}\).

\subsection{Meso-slow homogenization}
\label{secmsh}

To analyse the heterogeneous \pde{}s~\eqref{eq:shdifpde,eq:hdifpde} as a three-scale system we first establish and create the `meso-slow'  manifold homogenization.
That is, we use that the domain~\cD, \cref{fig:pattensemble}, is thin in~\(x_2\), that is \(\ell_2\ll\ell_1,L\).
Consequently, we consider the embedding \pde~\eqref{eq:emdifpde} by assuming the spatial variations in both~\(x_0\) and~\(x_1\) are on a relatively large length scale.  
That is, the derivatives of~\fu\ and~\(\kappa\) in both~\(x_0\) and~\(x_1\) are treated as small.

This homogenization is an example of that of a functionally graded material \cite[e.g.,][]{Anthoine2010} because the mesoscale variations in heterogeneity in~\(x_1\) are on a length scale \(\ell_1\gg\ell_2\)\,.  
Our slow manifold approach systematically incorporates all such variations within the analysis: for example, giving the fourth-order homogenization~\eqref{eqmshU}--\eqref{eqmsh}.
At orders higher than the leading order, we systematically find higher-order derivatives of the larger scale material structure: for example, the fourth-order coefficient in the upcoming \pde~\eqref{eqmsh} involves the third derivative~\(\fK'''(x_1)\) of the mesoscale variations.
This approach automatically incorporates effects obtained by the ``second-order homogenization'' \text{of \cite{Anthoine2010}.}

Akin to \cref{secsmh}, we here establish the basis of a slow manifold. Consider \pde~\eqref{eq:emdifpde}, with \(\D{x_0}{}\) and~\(\D{x_1}{}\) neglected, and in the Sobolov space~\(\HH_{[0,\ell_2]}\) with \(\ell_2\)-periodic boundary conditions in~\({x_2}\):
\begin{equation}
\D t\fu=\D {x_2}{}\left\{\kappa(x_1,x_2)\D {x_2}\fu \right\}.
\label{eq:emdifpde1}
\end{equation}
The basis established here applies at each and every~\(x_0,x_1\).

\paragraph{Equilibria}
For \pde~\eqref{eq:emdifpde1}, and for every~\(x_0,x_1\), in~\(\HH_{[0,\ell_2]}\) clearly \(\fu(t,x_2)=\fU\), constant in~\(x_2\), forms a subspace of equilibria~\EE.

\paragraph{Spectrum identifies invariant manifolds}
Since \pde~\eqref{eq:emdifpde1} is linear, the perturbation problem is identical at every equilibria in~\EE, namely~\eqref{eq:emdifpde1} itself.
We need to characterise the right-hand side operator in~\eqref{eq:emdifpde1}, \(\fL_0:=\D {x_2}{}\left\{\kappa(x_1,x_2)\D {x_2}\cdot \right\}\) with \(\ell_2\)-periodic boundary conditions, in~\(\HH_{[0,\ell_2]}\).

For example, in the case where the heterogeneity is constant with respect to~\(x_2\), a complete set of linearly independent eigenfunctions of~\(\fL_0\) are \(\e^{\i n k_2 x_2}\) for every integer~\(n\).
The spectrum of eigenvalues at~\((x_0,x_1)\) is then \(\lambda_{n}(x_1):=-\kappa(x_1)k_2^2n^2\).
This spectrum has one zero-eigenvalue and the rest are negative: \(\lambda_{n}(x_1)\leq -\kappa(x_1) k_2^2 <0\).
Hence, for every~\(x_0,x_1\), the non-zero eigenvalues satisfy \(\lambda_n(x_1)\leq-\kappa_{\min}k_2^2<0\).

For general heterogeneity, similar deductions to those of \cref{secsee}, establish that the operator~\(\fL_0\) is self-adjoint.
Hence the eigenvalues are real, and eigenfunctions orthogonal.
By a similar arguments to that of \cref{secsee}, it follows that the non-zero eigenvalues of the operator are bounded above by \(\lambda\leq-\beta_2\) for \(\beta_2:=\kappa_{\min}k_2^2=4\pi^2\kappa_{\min}/\ell_2^2\).

\paragraph{Slowly varying theory}
We proceed analogous to \cref{secrcvt} but now with two space dimensions (\(x_0\)~and~\(x_1\)) in which solutions are slowly varying.
\cite{Roberts2016a} developed rigorous slow manifold theory for spatio-temporal systems that have slow variations in multiple space dimensions. 
The operator~\(\fL_0\) satisfies the preconditions in Assumption~3 of \cite{Roberts2016a} with one zero eigenvalue (\(m=1,\alpha=0\)) and the rest negative.
By Proposition~1 of \cite{Roberts2016a}, and analogous to \cref{prosvpde}, we are assured that in a regime of slowly varying solutions, the mesoscale field~\(\fU(t,x_0,x_1)\) satisfies a \pde\ akin to~\eqref{eq:genhomopde}, but in~2D space, to a quantified error  \cite[eqn.~(52)]{Roberts2016a}, and upon neglecting exponentially decaying transients (decaying faster than~\(\e^{-\beta't}\) for \(\beta'\approx\beta_2\)).

A generating function argument \cite[\S\S3.2--3.3]{Roberts2016a}  establishes the validity of formal `slowly varying' analysis of the embedding \pde~\eqref{eq:emdifpde}.

\paragraph{An example mesoscale homogenization}
\cref{appmsh} lists and describes \cas\ code to construct the mesoscale homogenization for the example case of heterogeneous diffusivity
\begin{equation}
\kappa(x_1,x_2):=1/\big[1/\fK(x_1)+a_2\cos k_2x_2\big],
\quad\text{for }k_2:=2\pi/\ell_2\,,
\label{eqmshdiff}
\end{equation}
and some given mesoscale heterogeneity function~\(\fK(x_1)\).
The code derives that in terms of the mesoscale mean~\(\fU(t,x_0,x_1)\), and in terms of \(\fD_x:=\partial_{x_0}+\partial_{x_1}\), the microscale field (cf.~\eqref{eqqvgu}) is
\begin{align}
\fu&=\fU+\frac{\fK(x_1)a_2\sin k_2x_2}{k_2}\fD_x\fU
+\fD_x\left\{ \frac{\fK(x_1)a_2\cos k_2x_2}{k_2^2} \fD_x\fU \right\}
+\Ord{\fD_x^3}.
\label{eqmshU}
\end{align}
Simultaneously, the code derives that mesoscale mean field evolves according to the quasi-diffusion \pde
\begin{equation}
\D t\fU=\fD_x\left\{ \fK(x_1) \fD_x\fU \right\}
+\frac{a_2^2}{2k_2^2}\fD_x\left\{ \fK(x_1)^2 \fD_x^2\left[ \fK(x_1) \fD_x\fU \right] \right\}
+\Ord{\fD_x^6}.
\label{eqmsh}
\end{equation}
As to be expected, the effective diffusivity on the mesoscale,~\(\fK(x_1)\), is the harmonic mean over~\(x_2\) of the microscale diffusivity~\eqref{eqmshdiff}.
The theory of this subsection assures us that a truncation of the \pde~\eqref{eqmsh} models the dynamics of the embedded \pde~\eqref{eq:emdifpde} whenever and wherever the spatial gradients in~\(x_0\) and~\(x_1\) are `small' enough on the microscale~\(\ell_2\).
That is, \pde~\eqref{eqmsh} is \text{a mesoscale model.}

\subsection{Homogenization of the mesoscale model}

The mesoscale \pde~\eqref{eqmsh} is heterogeneous on the spatial length scale~\(\ell_1\) of variations in the coefficient~\(\fK(x_1)\).
For the case \(\ell_1\ll L\) we here derive the macroscale homogenization of this mesoscale model.
We have to choose some truncation of \pde~\eqref{eqmsh} to analyse further, so let's choose the leading order truncation \(\fU_t=\fD_x\left\{ \fK(x_1) \fD_x\fU \right\}\) as the mesoscale model.

To compare with the one-step homogenization, let's consider \pde~\eqref{eqmsh} in the case where the mesoscale diffusivity is
\begin{equation}
\fK(x_1):=1/\left[ 1+a_1\cos k_1x_1 \right]
\quad\text{for }k_1:=2\pi/\ell_1\,.
\label{eqhsmb}
\end{equation}

The theory and construction of the homogenization follows the same procedure as detailed before, so let's just summarise here.
First, we treat derivatives~\(\D{x_0}{}\) as `small'.  
Then the cross-sectional \pde\ is
\begin{equation}
\D t\fU=\D{x_1}{}\left\{ \fK(x_1) \D{x_1}\fU \right\},
\label{eqmshx}
\end{equation}
to be solved with \(\ell_1\)-periodic boundary conditions in~\(x_1\), and in space~\(\HH_{[0,\ell_1]}\).

Second, a subspace of equilibria of~\eqref{eqmshx} are \(\fU(t,x_1)=U\), constant in~\(x_1\).
The spectrum of the linear operator \(\cL_0:=\D{x_1}{}\left\{ \fK(x_1) \D{x_1}\cdot \right\}\) can be shown to be self-adjoint with one zero eigenvalue and the rest negative and bounded away from zero, \(\lambda\leq-\beta_1\), for \(\beta_1:=k_1^2\fK_{\min}=4\pi^2\fK_{\min}/\ell_1^2\).

Third, slowly varying theory for cylindrical domains, akin to \cref{secrcvt}, applies to assure us there exists a slow manifold homogenization of the mesoscale model in terms of a macroscale mean field~\(U(t,x_0)\).
Further, the slow manifold emerges exponentially quickly on a time-scale of~\(1/\beta_1\), and we may approximate the slow manifold via formal `slowly varying' analysis of the mesoscale \pde~\eqref{eqmsh}.

Fourth, the construction from the mesoscale \pde\  \(\fU_t=\fD_x\left\{ \fK(x_1) \fD_x\fU \right\}\) to the macroscale homogenization is the specific case of \cref{secxchm} with \(a_2=0\) and no variable~\(x_2\).
Consequently, corresponding results follow immediately, such as the slow manifold is~\eqref{eqqvgu} with \(a_2=0\), and the macroscale homogenization is \(U_t=U_{xx}\) to errors~\Ord{\partial_x^4}.

Here there is no point proceeding to higher order error in~\(\partial_x\) because we chose the leading order mesoscale \pde\ at the start of this subsection.
To obtain correct higher-order macroscale homogenizations, we would chose a correspondingly higher-order mesoscale homogenization.

\section{Conclusion}
\label{secconc}

This article further develops a novel theory and practice for mathematical\slash asymptotic homogenization to illuminate and resolve many issues in homogenization heterogeneous \pde{}s such as~\eqref{eq:hdifpde}.
\cref{secpse} invokes the technique of analysing an ensemble of all phase shifts, and extends it to multi-periodic cases, including the case of quasi-periodic.  
Future research may be able to extend the approach to interesting classes of random heterogeneity.
The key to this approach is to be able to recast the ensemble as a system which is mathematically homogeneous \text{in the macroscale.}

Given a now mathematically homogeneous macroscale, \cref{secsmh} applies recent developments of dynamical systems theory to newly establish the homogenization as a rigorous slow manifold of the given \pde\ system.
No~``\(\epsilon\)'' is required so the results hold for \emph{finite} scale separation.
No variational principles are required, so the approach applies to a wide range of physical problems, both dissipative and wave-like.
The dynamical systems approach provides a framework that future research and applications may exploit to rigorously derive accurate initial conditions, boundary conditions, forcing projections, uncertainty evaluation, \text{and error bounds.}

This dynamical systems framework provides a very practical way to construct homogenized models, as in the example of \cref{secxchm,seccachm}.
It shows the ease in using computer algebra to perform and check the homogenization~\eqref{eqsqvg}.
The method systematically and straightforwardly extends to higher-order to derive so-called `second-order homogenization': indeed, \cref{seccaisc} simply proceeds further to explicitly construct the `fourth-order homogenization'~\eqref{eqimphpde}.

One advantage in the computer algebra is that, for some classes of heterogeneities, we can proceed to very high-order in the analysis.
\cref{secsrh} computed to 34th-order and uses the resulting generalised homogenization to quantify how `much smaller' the microscale really has to be compared to the macroscale.
Remarkably, the evidence is that, ideally, macroscale lengths can be resolved down to just twice the microscale length! 
Of course, in practice, there are many confounding aspects, but this factor of two is the mathematically \text{ideal bound.}

The spatio-temporal resolution of homogenizations could be improved, in future research, by multi-modal\slash multi-zonal homogenization of the heterogeneous microscale system by adapting approaches developed for shear dispersion \cite[]{Strunin01a, Watt94b}.

Lastly, \cref{sec3sh} explores how the approach developed here includes the case when there are three separated length scales in the physical problem: a microscale, mesoscale, and macroscale.
The exploration shows how one can analyse such a system in one step from the micro+mesoscale to the macroscale, or one can analyse the system in two steps from micro to meso, and then from meso to macro.
Because sound modelling is transitive, the resulting homogenizations are the same. 
Further, if one's ultimate aim is a macroscale spatial discretisation in order to compute predictions, that is the aim is not really the homogenized \pde, then instead of the two step process of firstly homogenization and secondly spatial discretisation, future research could do it all in one step from microscale heterogeneous \pde\ to macroscale discretisation using a process like that introduced in a shear dispersion example \cite[]{MacKenzie03}.

\paragraph{Acknowledgements}
I thank Arthur Norman and colleagues who maintain the computer algebra software Reduce.  
The Australian Research Council Discovery Project grants DP200103097 and DP220103156 helped support \text{this research.}

\appendix
\section{Computer algebra to construct an homogenized PDE} 
\label{seccachm}

The following code constructs the
homogenization~\eqref{eq:genhomopde}, discussed by
\cref{secxchm}, of the
\pde{}s~\eqref{eq:hdifpde,eq:emdifpde} for the family of
problems with diffusivity~\eqref{egkappa} for microscale
periodicities~\(\ell_1,\ell_2\).  Let's use the computer
algebra system Reduce\footnote{%
\url{http://www.reduce-algebra.com}} as it is free,
flexible, and fast~\cite[e.g.]{Fateman2002}.  First improve
printed output:
\begin{reduce}
on div; off allfac; on revpri;
factor d,df;
\end{reduce}
Let the arguments of trig function heterogeneity be denoted 
by~\(q_i\) and define wavenumbers \(k_i:=2\pi/\ell_i\), 
so the diffusivity~\(\kappa(x_1,x_2)\) is the following.
\begin{reduce}
depend q1,x1,z; depend q2,x2,z;
let { df(q1,x1)=>k1, df(q2,x2)=>k2 };
kappa:=1/(1+a1*cos(q1)+a2*cos(q2));
\end{reduce}

Write approximations to the slow manifold model of the
embedding \pde~\eqref{eq:emdifpde} in terms of a `mean'
field~\(U(t,x)\), denoted by~\verb|uu|, that evolves
according to \(\D tU=\verb|dudt|\) for whatever \verb|dudt|
happens to be.
\begin{reduce}
depend uu,x,t;
let df(uu,t)=>dudt;
\end{reduce}
This code uses~\(x\) for~\(x_0\) as the two are synonymous
in practice.

\subsection{Quickly verify a leading approximation}
\label{seccaqvg}
First, we guess and then verify the approximate
field~\eqref{eqqvgu}. This is coded, with ordering parameter
\(\verb|d|\) that counts the number of \(x\)-derivatives
in~\(U\), as
\begin{reduce}
u:=uu+d*(a1/k1*sin(q1)+a2/k2*sin(q2))*df(uu,x)
     +d^2*(a1/k1^2*cos(q1)+a2/k2^2*cos(q2))*df(uu,x,2);
\end{reduce}
such that~\(U\) evolves according to the homogenized \pde\
\(\D tU = \DD xU +\cdots\) (the coefficient~\(1\) of the
mean diffusivity is the classic harmonic mean of the
specific heterogeneous diffusivity~\eqref{egkappa}):
\begin{reduce}
dudt:=d^2*df(uu,x,2);
\end{reduce}
To verify, substitute into the governing
\pde~\eqref{eq:emdifpde} and find the \pde's residual is
zero to an error~\Ord{\partial_x^3}---counting the
derivatives with the order parameter~\verb|d|.
\begin{reduce}
let d^3=>0; 
flux:=-kappa*(d*df(u,x)+df(u,x1)+df(u,x2))$
pde:=df(u,t)+d*df(flux,x)+df(flux,x1)+df(flux,x2);
\end{reduce}
Since the \pde\ residual is~\Ord{\partial_x^3}, then the
slow manifold~\eqref{eqqvgu} is correct to
error~\Ord{\partial_x^3} \cite[]{Roberts2013a}.

The next subsection repeatedly computes the residual, to
drive corrections, until the residual is zero to a specified
order of error, and hence gives a slow manifold to the same
order of error.

\subsection{Iteration systematically constructs}
\label{seccaisc}

Second, we iteratively construct the improved
homogenization~\eqref{eqimphpde}. This code executes in
about four~seconds.
\begin{reduce}
write "
Second, Iteratively Construct
-----------------------------";
\end{reduce}
The terms~\Ord{\partial_x^4} appear, at first, impossible to
find with exact algebra, so instead we construct in a power
series in the amplitude of the heterogeneity.  Let~\verb|a|
count the number of amplitude factors arising
from~\(\kappa\).  Choose the following order of errors in
both derivatives and \(a:=\sqrt{a_1^2+a_2^2}\). It
eventuates that for this particular
heterogeneity~\(\kappa\), the results to
error~\Ord{\partial_x^{\texttt{ord}}} appear exact in
heterogeneity amplitude~\(a\).
\begin{reduce}
ord:=7;
\end{reduce}
Here assume heterogeneity~\(\kappa\) is a reciprocal, but
could be more general.
\begin{reduce}
rkappa:=sub({a1=a*a1,a2=a*a2},1/kappa);
\end{reduce}

\paragraph{Optional} We may set \(a_2=0\) to get simpler
results for a single periodic component, and change
truncation to the same order, now high-order.  Order~35 takes
roughly one minute compute time. The high-order
coefficients of the homogenized \pde\ resulting from this
option are exact, as far as they go, because the evidence is
that the expansions in homogeneity amplitude~\(a\) extend to
at most the same order as the order of derivatives.
\begin{reduce}
if 1 then begin a2:=0; a1:=k1:=1; ord:=35; end;
\end{reduce}

\paragraph{In general}
Construct heterogeneous diffusivity~\(\kappa\) as the
reciprocal of~\verb|rkappa|, presuming \(\kappa\approx 1\):
\begin{reduce}
kappa:=1+for n:=1:ord-1 sum (1-rkappa)^n$
res:=(kappa*rkappa-1 where a^~p=>0 when p>=ord);
if res neq 0 then rederr("kappa reciprocal error");
\end{reduce}

Start the iteration from the trivial approximation for the
field and its evolution.
\begin{reduce}
u:=uu;
dudt:=0;
\end{reduce}
Seek solution to the specified orders of errors.
\begin{reduce}
for it:=1:999 do begin write "
**** ITERATION ",it;
\end{reduce}

Progressively truncate the order of the order parameter so
that we control the residuals better: the bound in this
if-statement is the aimed for ultimate order of error.  
\begin{reduce}
  if it<ord then let {d^(it+1)=>0, a^(it+1)=>0};
\end{reduce}
Compute the \pde\ residual via the flux, trace printing the
length of the residual expression:
\begin{reduce}
  flux:=-kappa*(d*df(u,x)+df(u,x1)+df(u,x2));
  pde:=df(u,t)+d*df(flux,x)+df(flux,x1)+df(flux,x2);
  pde:=trigsimp(pde,combine);
  write lengthpde:=length(pde);
\end{reduce}
Update the evolution via the solvability condition from
the mean over~\(x_1,x_2\):
\begin{reduce}
  gd:=-(pde where {sin(~a)=>0,cos(~a)=>0});
  if it<6 then write gd:=gd else write lengthgd:=length(gd);
  dudt:=dudt+gd;
  rhs:=pde+gd;
\end{reduce}

Attempt to solve the update in two steps, each step with two
integrals. Use the following operator, given we already
made~\verb|qi| depend upon dummy~\verb|z|.
\begin{reduce}
  if it=1 then begin
    operator intx; linear intx;
    let { intx(cos(~q),z) => sin(q)/(df(q,x1)+df(q,x2))
        , intx(sin(~q),z) =>-cos(q)/(df(q,x1)+df(q,x2))
        };
  end;
\end{reduce}
\begin{itemize}
\item First integrate \(v_{x_1}+v_{x_2}=\verb|rhs|\) with
zero mean. If any \verb|intx(1,z)| appear then mistake.
\begin{reduce}
  v:=trigsimp( intx(rhs,z) );
  if not freeof(v,intx(1,z)) then rederr("ABORT");
\end{reduce}
\item Second, solve \(\kappa u'_z =\kappa (u'_{x_1}
+u'_{x_2}) =v\), using divide by~\(\kappa\) and integrate
again. Assume any presence of \verb|intx(1,z)| is an error
that eventually sorts itself out in the iteration, so
obliterate here! The update then has mean zero.
\begin{reduce}
  ud:=intx(trigsimp(v*rkappa,combine),z);
  ud:=(ud where intx(1,z)=>0);
\end{reduce}
\end{itemize}
Add in the update.  Perhaps trigsimp is useful here.
\begin{reduce}
  u:=u+trigsimp(ud);
\end{reduce}

Exit iterative for-loop when the \pde\ residual is zero
\begin{reduce}
  if pde=0 then write "Success: ", it:=it+10000;
end;
showtime;
if pde neq 0 then rederr("Iteration failure");
\end{reduce}

Write out homogenized \pde, optionally its coefficients, and
end the script.
\begin{reduce}
dudt:=dudt;
if ord>7 then begin cs:=coeff(dudt,d);
    cs:=for n:=3 step 2 until length(cs) 
        collect part(cs,n)/df(uu,x,n-1);
    on rounded; off nat; write cs:=cs;  
    on nat; off rounded; end;
end;
\end{reduce}


\section{Meso-scale homogenization code}
\label{appmsh}

Construct the asymptotic expansion of the meso-slow
homogenization of three-scale diffusion supported by the
theory of \cref{secmsh}.
\begin{reduce}
on div; off allfac; on revpri;
factor d,b,a2,k2;
\end{reduce}
The diffusivity~\(\kappa(x_1,x_2)\) is the following, where
\verb|b(n)| denotes the \(n\)th derivative of~\(b(x_1)\).
At the end we recast in \(\verb|c|=\fK(x_1):=1/b(x_1)\).
\begin{reduce}
depend kappa,x1,x2;
operator b; depend b,x1; 
kappa:=1/(b(0)+a2*cos(k2*x2));
let { df(b(~n),x1) => b(n+1) };
\end{reduce}

Write approximations to the slow manifold model of the
embedding \pde~\eqref{eq:emdifpde} in terms of a `mean'
field~\(U(t,x)\) that evolves according to \(\D
tU=\verb|dudt|\).
\begin{reduce}
depend uu,x,x1,t;
let df(uu,t)=>dudt;
\end{reduce}
This code uses~\(x\) for~\(x_0\) as the two are synonymous
in practice. Start from the leading approximation for the
field and its evolution.
\begin{reduce}
u:=uu$
dudt:=0$
\end{reduce}

Iteratively construct to this order of error in slow
derivatives.
\begin{reduce}
ordd:=5;
for it:=1:999 do begin write "
**** Iteration ",it;
\end{reduce}
Progressively truncate the order of the order parameter so
that we better control the residuals.  
\begin{reduce}
  if it<ordd then let d^(it+1)=>0;
\end{reduce}
Compute the \pde\ residual via the flux, trace printing the
length of the residual expression:
\begin{reduce}
  flux:=-kappa*(d*df(u,x)+d*df(u,x1)+df(u,x2));
  pde:=df(u,t)+d*df(flux,x)+d*df(flux,x1)+df(flux,x2);
  pde:=trigsimp(pde);
  write lengthpde:=length(pde);
\end{reduce}

Solvability updates the evolution:
\begin{reduce}
  dudt:=dudt+(gd:=-int(pde,x2,0,2*pi/k2)/(2*pi/k2));
\end{reduce}
Update the field via some integrals:
\begin{reduce}
  ud:=trigsimp(int(int(pde+gd,x2)/kappa,x2));
  udx:=sub(x2=2*pi/k2,ud)-sub(x2=0,ud); 
  u:=u+ud-int(ud,x2,0,2*pi/k2)/(2*pi/k2)
         -udx*(x2-pi/k2+a2/k2/b(0)*sin(k2*x2))/(2*pi/k2);
\end{reduce}

Exit iterative for-loop when the residual is zero
\begin{reduce}
  showtime;
  if pde=0 then write "Success: ", it:=it+10000;
end;
if pde neq 0 then rederr("Iteration failure");
\end{reduce}

Write out homogenized \pde.
\begin{reduce}
dudt:=dudt;
\end{reduce}

Recast in terms of \(\verb|c|=\fK(x_1):=1/b(x_1)\).
\begin{reduce}
operator c; depend c,x1; factor c;
let { df(c(~n),x1)=>c(n+1), b(~n)=>df(1/c(0),x1,n) };
dudt:=dudt;
\end{reduce}
Check compact form of the evolution \pde, and finish the script.
\begin{reduce}
procedure dx(a); d*(df(a,x)+df(a,x1))$
err:=-(dudt where d^5=>0) +dx(c(0)*dx(uu))
    +a2^2/k2^2/2*dx(c(0)^2*dx(dx(c(0)*dx(uu))));
if err neq 0 then rederr("simplification failure");
end;
\end{reduce}


\end{document}